\documentclass[12pt]{amsart}

\usepackage[dvips]{graphicx}
\usepackage{amsfonts,amssymb,amscd,bm}
\newcommand{\R}{{\mathbb R}}
\newcommand{\Z}{{\mathbb Z}}
\newcommand{\ph}{{\varphi}}
\newcommand{\ep}{{\varepsilon}}

\newcommand{\al}{{\alpha}}

\newcommand{\La}{{\Lambda}}

\newcommand{\Si}{{\Sigma}}

\newcommand{\ti}{\tilde}
\newcommand{\id}{{\mathrm{id}}}
\newcommand{\pr}{{\mathrm{pr}}}
\newcommand{\tg}{{\mathrm{TG}}}
\newcommand{\ct}{{\mathrm{CT}}}
\newcommand{\up}{{\mathrm{UP}}}
\newcommand{\tb}{{\mathrm{TB}}}
\newcommand{\hb}{{\mathrm{HB}}}
\newcommand{\SK}{{\mathrm{SK}}}
\newcommand{\SL}{{\mathrm{SL}}}
\newcommand{\CS}{{\mathrm{CS}}}

\newcommand{\bd}{{\partial}}
\newcommand{\bu}{{\bullet}}

\theoremstyle{definition}

\title[2-dimensional links live in a universal 3-dimensional polyhedron]
{All 2-dimensional links in 4-space live inside\\ 
 a universal 3-dimensional polyhedron}

\author[kearton]{C.~Kearton}
\address{Department of Mathematical Sciences,   
 Durham University, Durham DH1 3LE, United Kingdom}
\email{ cherry.kearton@durham.ac.uk }

\author[kurlin]{V.~Kurlin}
\address{Department of Mathematical Sciences,   
 Durham University, Durham DH1 3LE, United Kingdom}
\email{ vitaliy.kurlin@durham.ac.uk }

\date{First version: January 23, 2008; this version: April 7, 2008. }

\begin{document}

\maketitle

\begin{abstract}
The hexabasic book is the cone of the 1-dimensional skeleton
 of the union of two tetrahedra glued along a common face.
The universal 3-dimensional polyhedron $\up$
 is the product of a segment and the hexabasic book.
We show that any 2-dimensional link in 4-space is isotopic
 to a surface in $\up$.
The proof is based on a representation of surfaces in 4-space 
 by marked graphs, links with double intersections in 3-space.
We construct a finitely presented semigroup whose central elements
 uniquely encode all isotopy classes of 2-dimensional links. 
\end{abstract}

%===================================================

\section{Introduction}%1

%--------------------------------------------------------------------------------------

\subsection{Brief summary}%1.1
\noindent
\smallskip

This is a research on the interface between 
 geometric topology, singularity theory and semigroups.
A 2-link is a closed 2-dimensional surface in 4-dimensional space $\R^4$.
We study 2-links up to isotopy that is a smooth deformation of 
 the ambient 4-dimensional space.
We prove that any 2-link is isotopic to a surface embedded
 into the universal 3-dimensional polyhedron $\up$.
We also reduce the isotopy classification of 2-links in 4-space
 to a word problem in a finitely presented semigroup.
\medskip

%--------------------------------------------------------------------

\subsection{The universal polyhedron containing 2-dimensional links}%1.2
\noindent
\smallskip

First we define the universal 3-dimensional polyhedron $\up$.
\medskip

\noindent
{\bf Definition 1.1.}
The \emph{theta} graph $\tg$ consists of 3 edges connecting 2 vertices.
The \emph{circled} theta graph $\ct$ is $\tg\cup S^1$, where
 the circle $S^1$ meets each edge of $\tg$ in one point, see Fig.~1.
Then $\ct$ is the 1-dimensional skeleton of two tetrahedra glued
 along a common face.
The \emph{hexabasic book} $\hb$ is the cone of $\ct$.
Being embedded in 3-space, the book $\hb$ divides 
 a neighbourhood of the central vertex into 6 parts.
The \emph{universal} 3-dimensional polyhedron
 is $\up=\hb\times[-1,1]$.

\begin{figure}[!h]
\includegraphics[scale=1.0]{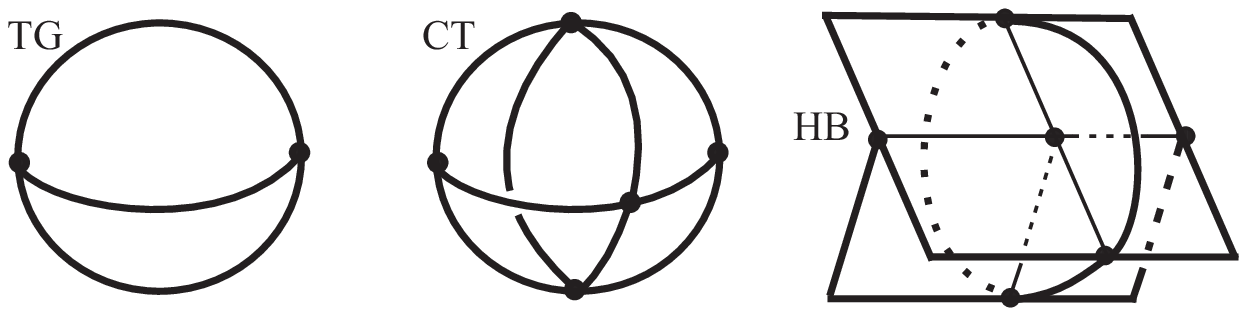}
\caption{The theta graph $\tg$, circled theta graph $\ct$, book $\hb$}
\end{figure}

We will work in the smooth category, i.e. all diffeomorphisms are $C^{\infty}$-smooth.
We will make necessary comments on similar constructions in the PL case.
\medskip
 
\noindent
{\bf Definition 1.2.}
An \emph{embedding} is a diffeomorphism onto its image.
A \emph{2-link} is a closed (possibly disconnected or 
 non-orientable) smooth surface $S$ embedded into $\R^4$.
%If the surface $S$ is connected then the image 
% $S\subset\R^4$ is called a \emph{2-knot}.                           
An \emph{isotopy} between 2-links $S$ and $S'$ is 
 a continuous family of diffeomorphisms $F^u:\R^4\to\R^4$, 
 $u\in[0,1]$, such that $F^0=\id_{\R^4}$, $F^1(S)=S'$.
\medskip
 
%\noindent
Fix the 4th coordinate $t$ in 4-space $\R^3\times \R$.
Then a 2-link in $\R^3\times\R$ can be studied in terms of 
 its \emph{cross-sections} $S_t=S\cap(\R^3\times\{t\})$ \cite{FM}.
Any 2-link can be isotopically deformed to
 a surface $S\subset\R^3\times[-1,1]$ such that the projection
 $\pr:S\to[-1,1]$ has distinct non-degenerate critical values.
A general cross-section $S_t$ is a classical link in $\R^3\times\{t\}$, while 
 a cross-section containing a saddle is a link with a double point.
When $t$ passes through a saddle,
 the cross-section $S_t=S\cap(\R^3\times\{t\})$ changes  
 by the Morse modification in the left picture of Fig.~2.
\smallskip

\begin{figure}[!h]
\includegraphics[scale=1.0]{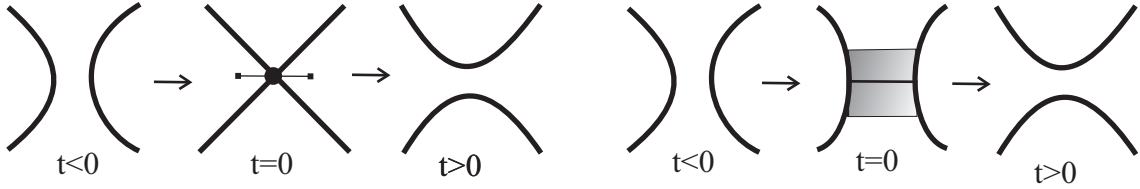}
\caption{Resolving a singular point and a band in $\R^3$}
\end{figure}

A PL analogue of the smooth approach is to decompose 
 a 2-link $S\subset\R^3\times[-1,1]$ into handles located
 in different sections $\R^3\times\{t_j\}$.
The 1-handles of $S$ will be represented by bands that have 
 a distinguished core and are attached to a classical link in 3-space.
Any attached band can be retracted to a singular point 
 marked by a bridge encoding the core of the band.
The cross-sections of $S$ below and above 
 every 1-handle locally look like the right picture of Fig.~2.
\noindent
\medskip
 
%------------------------------------------------------------------------------------------------

\subsection{Main results}%1.3
\noindent
\smallskip

\noindent
The hexabasic book $\hb$ is closely related to the 
 \emph{3-page book} $\tb$, the cone of the \emph{theta} graph $\tg$
 consisting of 3 edges connecting 2 vertices, see Fig.~7.
The \emph{binding} segment of $\tb$ is the cone of the 2 vertices of $\tg$.
From another point of view, the 3-page book $\tb$ can be considered as $\R\times T$,
 where $T$ is the \emph{triod} consisting of 3 edges connecting 
 the central vertex $O$ to other 3 vertices, here the binding axis $\al$ is $\R\times O$.
The hexabasic book $\hb$ is obtained from $\tb$ by adding
 3 half-disks whose 6 boundary radii are attached to 
 the 3 edges of $\{0\}\times T$, see Fig.~1.
\medskip

\noindent
{\bf Theorem~1.3}.
Any 2-dimensional link $S\subset\R^4$ is isotopic to a surface embedded 
 into the universal 3-dimensional polyhedron $\up=\hb\times[-1,1]$.
\medskip

%\noindent
The key idea of Theorem~1.3 is to put a given $S$ surface in general position
 and consider its cross-sections $S_t$ through saddles of $\pr:S\to[-1,1]$,
 see Claim~2.3.
Such a cross-section $S_t$ is a link with exactly one singular point, so
 $S_t$ can be embedded into the 3-page book $\tb$ using 
 the technique of 3-page embeddings developed in \cite{Kur,KV},
 see Proposition~3.2.
Both resolutions of the singular point of $S_t$ can be realised in $\tb$, 
 i.e. the embedding extends to a regular neighbourhood of $S_t$ in $S$.
It remains to embed the complement of the regular neighbourhoods of all saddles
 into $\hb\times[-1,1]$ realising any isotopy of classical links in $\hb$, 
 see Lemma~3.4.
\medskip

We will develop a 1-dimensional calculus for 2-links as follows.
Any 2-link $S$ in general position in $\R^3\times[-1,1]$ 
 can be represented by a banded link $BL$ whose bands are 
 associated to the saddles of $\pr:S\to[-1,1]$, see Proposition~2.6(i).
Retracting each band to a point, we get a marked graph whose
 singular points are marked by bridges encoding the cores of bands.
There is a complete set of moves on marked graphs generating
 any isotopy of 2-links in 4-space, see Proposition~4.2.
Any marked graph can be embedded into the 3-page book $\tb$
 and can be encoded by a word in the alphabet of 15 letters.
The moves on marked graphs are translated into relations on words,
 which leads to the universal semigroup $\SL$ of 2-links in 4-space.
\smallskip

%\noindent
Introduce the universal semigroup $\SL$ generated by
 the letters $a_i,b_i,c_i,d_i,x_i$ subject to relations (1)-(8), 
 where  $i\in\Z_3=\{0,1,2\}$, e.g. $0-1=2\pmod{3}$.
\medskip

\noindent
(1) $d_0d_1d_2=1$, $\quad b_id_i=1=d_ib_i$;
\smallskip

\noindent
(2) $a_i=a_{i+1}d_{i-1}$, $\quad b_i=a_{i-1}c_{i+1}$, 
     $\quad c_i=b_{i-1}c_{i+1}$, $\quad d_i=a_{i+1}c_{i-1}$;
\smallskip

\noindent
(3) $uv=vu$, $u\in\{a_ib_i, d_ic_i, b_{i-1}d_id_{i-1}b_i, d_ix_ib_i\}$,
     $v\in\{a_{i+1}, b_{i+1}, c_{i+1}, b_id_{i+1}d_i, x_{i+1}\}$;
\smallskip

\noindent
(4) $x_{i-1}=b_{i+1}x_id_{i+1}$, 
     $\quad b_i x_i b_i=a_i (b_i x_i b_i) c_i$,
     $\quad d_i x_i d_i=a_i (d_i x_i d_i) c_i$;
\smallskip

\noindent
(5) $(d_i x_i b_i)d_i^2d_{i+1}^2d_{i-1}^2=d_i^2d_{i+1}^2d_{i-1}^2(d_i x_i b_i)$;
\smallskip

\noindent
(6) $a_i x_i=a_i$, $a_i b_i x_i d_i c_i =1$;
\smallskip

\noindent
(7) $d_i x_i b_i c_i x_i=b_i x_i d_i c_i x_i$;
\smallskip

\noindent
(8) 
$w_i d_{i+1} d_i^2 d_{i-1} a_{i+1} b_{i+1} x_i b_i d_{i+1} b_i^2 b_{i+1} d_i^2=
 w_i b_{i-1} b_i a_i b_{i+1} a_{i+1} d_i^2 c_{i-1} b_i x_i b_i$,\\
 where $w_i=a_i b_i x_i b_i c_i$. 
%(a_ib_ix_ib_ic_ib_ix_i)d_ib_{i-1}^2c_i^2d_{i-1}^2=
%       (a_ib_ix_ib_ic_ib_ix_i)d_{i+1}b_{i-1}d_i^2c_{i-1}^2b_i^2$.
\bigskip

One of the 6 relations $b_id_i=1=d_ib_i$ is superfluous and can be deduced
 from the remaining relations in (1).
Moreover, the commutativity of $d_ic_i$ with $a_{i+1},b_{i+1}$
 follows from the other relations in (3), see more details in \cite{Kur}.
So the semigroup $\SL$ is generated by 15 letters and 96 relations.
\bigskip

\noindent
{\bf Theorem~1.4}.
Any 2-link $S\subset\R^4$ is encoded by an element
 $w_S\in\SL$ in such a way that 2-links $S,S'$ are isotopic
 if and only if their encoding elements $w_S$ and $w_{S'}$
 are equal in $\SL$.
An element $w\in\SL$ encodes a 2-link if and only if $w$ is central in $\SL$.
\medskip

\noindent
\emph{Outline.}
In section~2 one represents 2-links in 4-space 
 by banded links and marked graphs in 3-space.
Theorems~1.3 and 1.4 are proved in sections~3 and 4, respectively.
Banded links are more convenient for deriving a complete set of moves
 generating any isotopy of 2-links.
Marked graphs will be used to prove our main results on embedding
 and encoding 2-links up to isotopy.
\smallskip

\noindent
{\bf Acknowledgement.}
The authors thank S.~Carter, F.~Tari for useful discussions.
\medskip

%===================================================

\section{Representing 2-links by banded links and marked graphs}%2

%------------------------------------------------------------------------------------------------

\subsection{Critical level embeddings of 2-links in 4-space}%2.1
\noindent
\smallskip

\noindent
Here we describe the PL approach where a 2-link is isotopically deformed
 to a nice embedding with handles at different levels.
The smooth version of crucial Claim~2.3(ii) is a standard statement 
 on general position proved in Appendix.
\medskip

\noindent
{\bf Definition 2.1.}
A \emph{handle} of dimension $n$ and index $k$ is $D^k\times D^{n-k}$.
% the attaching area is $\bd H=\bd D^k\times D^{n-k}$.
A \emph{handle decomposition} of a manifold $M^n$ is a sequence
 of submanifolds $M_0\subset M_1\subset\dots\subset M_l=M$, where
 $M_0$ is a disjoint union of $n$-dimensional disks, 
 each $M_{i+1}$ is obtained from $M_i$
 by adding a handle of some index $k_i$.
Formally, one can express $M_{i+1}=M_i\cup_{\ph_i}(D^{k_i}\times D^{n-k_i})$, 
 where $\ph_i:\bd D^{k_i}\times D^{n-k_i}\to\bd M_i$ is an embedding.
If before and after each handle addition one inserts a \emph{collar}, 
 the product of the attaching area and a segment, then
 one gets a \emph{collared handle} decomposition \cite[p.~416]{KL}. 
\medskip

\begin{figure}[!h]
\includegraphics[scale=1.0]{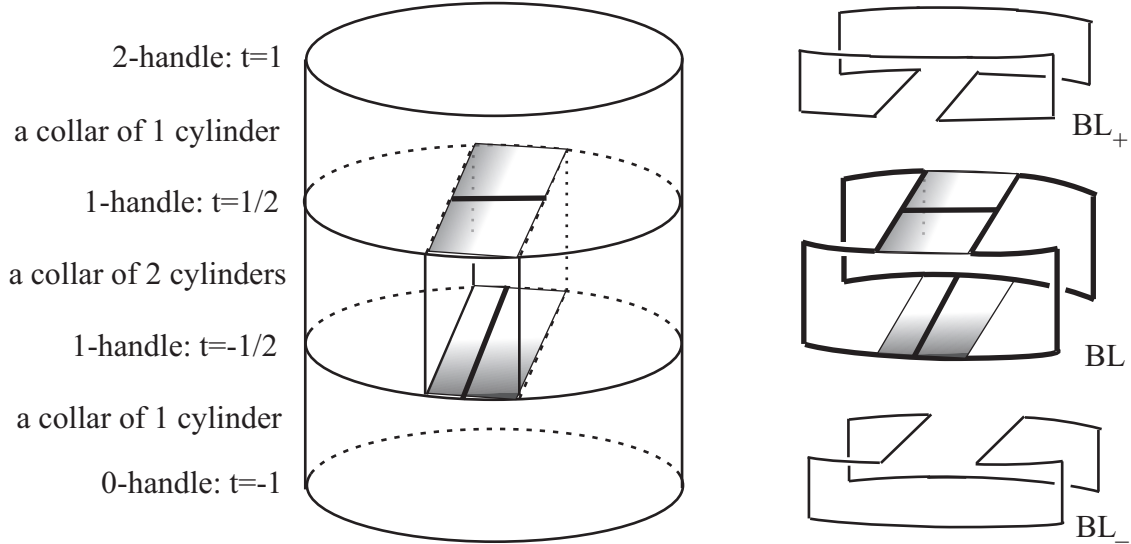}
\caption{A critical level embedded torus and its banded link $BL\subset\R^3$}
\end{figure}
\smallskip

\noindent
A 2-link with a collared handle decomposition
 can be nicely embedded in $\R^4$.
The left picture of Fig.~3 shows a similar embedding, 
 where the standard 2-torus in $\R^3$ 
 has the collared handle decomposition 
 consisting of 4 handles and 3 collars:
\smallskip

\noindent
1) the lowest handle is a 0-handle (a disk) at the level $t=-1$;
\smallskip

\noindent
2) the 2 intermediate handles are 1-handles (bands) at the levels $t=\pm 1/2$;
\smallskip

\noindent
3) the highest handle is a 2-handle (a disk) at the level $t=+1$.
\bigskip

\noindent
{\bf Definition 2.2.}
A \emph{critical level} PL embedding is a PL embedding of 
 a 2-link $S\subset\R^3\times[-1,1]$ with 
 a collared handle decomposition satisfying (i), (ii) \cite[p.~417]{KL}: 
\smallskip

\noindent
(i) the handles are in different sections $\R^3\times\{t_j\}$, 
 where $-1<t_1<\dots<t_n<1$;
\smallskip

\noindent
(ii) each collar between adjacent handles of $S$ is embedded as 
 the direct product $A\times[t_j,t_{j+1}]\subset\R^3\times[t_j,t_{j+1}]$, 
 where $A\subset\R^3$ is the attaching area of the handles.
\medskip

\noindent
A smooth embedding $S\subset\R^3\times[-1,1]$ 
 is called a smooth \emph{critical level} embedding
 if the projection $\pr:S\to[-1,1]$ has all its critical points
 in different sections $\R^3\times\{t_j\}$.
This is a general position assumption.
\medskip

%The idea to use cross-sections for representing 2-links goes back to R.~Fox.

\noindent
{\bf Claim~2.3.} (i) \cite[Theorem~1, p.~420]{KL}
Any 2-dimensional PL link in 4-space is isotopic to the image of 
 a critical level PL embedding $S\subset\R^3\times[-1,1]$.
\smallskip

\noindent
(ii) Any smooth 2-link is smoothly isotopic to
 a surface $S\subset\R^3\times[-1,1]$ such that all critical points of 
 $\pr:S\to[-1,1]$ are non-degenerate and have distinct values.
\qed
\medskip

We will use the smooth version of Lemma~2.3(ii), which will be 
 deduced from the transversality theorem of Thom in Appendix.
Claim~2.3(i) is worth keeping in mind when one associates 
 a banded link to a 2-link in Proposition~2.6(i).
\smallskip

%------------------------------------------------------------------------------------------------

\subsection{Representing 2-links in 4-space by banded links in 3-space}%2.2
\noindent
\smallskip

\noindent
We define banded links, links with bands, which will represent 2-links in 4-space.
\medskip

\noindent
{\bf Definition 2.4.}
A \emph{banded} link is a collection of circles and bands in $\R^3$ such that
\smallskip

\noindent
(i) the circles and bands are non-oriented and non-self-intersecting;
\smallskip

\noindent
(ii) the circles and bands are disjoint except for each band having 
 a pair of opposite sides \emph{attached} to disjoint arcs in the circles,
 the other sides are called \emph{free}.
\smallskip

\noindent
In every band we mark its \emph{core}, 
 an arc connecting its attached opposite sides, see Fig.~3.
Banded links are considered up to isotopy of $\R^3$.
The bands of a banded link will represent 1-handles of a 2-link.
In every band $B$ of a banded link $BL$ consider the opposite free sides 
 not connected by the core of $B$.
Replace $B$ by its free sides, the resulting usual non-oriented link in $\R^3$
 is called the \emph{positive} resolution $BL_+$ of the banded link $BL$, see Fig.~3.
Similarly define the \emph{negative} resolution $BL_-$ replacing 
 every band $B$ by the opposite attached sides connected by the core of $B$.
A banded link $BL$ is \emph{admissible}, if 
 both resolutions $BL_{\pm}$ are trivial links.
\medskip

%\noindent
If a PL 2-link $S\subset\R^3\times[-1,1]$ 
 has all its 1-handles in the zero section $\R^3\times\{t=0\}$,
 then the cross-section $S_0=S\cap(\R^3\times\{t=0\})$ is a banded link.
We will use much weaker assumptions and construct
 a banded link for any critical level embedding.
Proposition~2.6 leads to a calculus for 2-links in Proposition~4.2 and 
 provides a function from the set of 2-links to the set of admissible banded links.
\medskip

\noindent
{\bf Definition 2.5.}
Given a 2-dimensional surface $S$, consider the space of all smooth functions 
 $f:S\to\R^4$ with the Whitney topology, see Appendix.
The space $\CS$ of all 2-links $S\subset\R^4$ has the induced topology.
Points in $\CS$ 
% are considered as 2-links and 
 will be classified using the projection $\pr:S\to\R$ to the 4th coordinate $t$.
%are also denoted by $S$.
A 2-link $S\in\CS$ is called
\smallskip

\noindent
$\bu$
\emph{generic} if all critical points of $\pr$ 
 are non-degenerate and have distinct values;
\smallskip

\noindent
$\bu$
an \emph{$A_1^+ A_1^+$-singularity} if $S$ fails to be generic because of 
 2 non-degenerate extrema of $\pr:S\to\R$ that have the same value;
\smallskip

\noindent
$\bu$
an \emph{$A_1^+ A_1^-$ -singularity} if $S$ fails to be generic because of 
 a non-degenerate saddle and extremum of $\pr:S\to\R$ that have the same value;
\smallskip

\noindent
$\bu$
an \emph{$A_1^- A_1^-$ -singularity} if $S$ fails to be generic because of 
 2 non-degenerate saddles of $\pr:S\to\R$ that have the same value;
\smallskip

\noindent
$\bu$
an \emph{$A_2$-singularity} if $S$ fails to be generic because of 
 a singularity of $\pr:S\to\R$ having the form $\pr(x,y)=x^2-y^3$
 in local coordinates $x,y$.
\medskip

The sign in the notation above is the sign of the determinant 
 $\pr_{xx}\pr_{yy}-\pr_{xy}^2$ of the Jacobi matrix of 2nd order derivatives
 at a critical point.
Denote by $\Si_{++},\Si_{+-},\Si_{--}$ and $\Si_{2}$ the \emph{subspaces} 
 of the corresponding singularities in the space $\CS$.
Introduce the \emph{singular} subspace $\Si=\Si_{++}\cup\Si_{+-}\cup\Si_{--}\cup\Si_2$.
An isotopy of 2-links can be considered as a path in $\CS$.
In Proposition~2.6 we consider paths 
 nicely meeting the singular subspace $\Si$.
\medskip

\noindent
{\bf Proposition 2.6.}
(i) To any a critical level embedding $S\subset\R^3\times[-1,1]$
 we associate a banded link $BL$ well-defined up to 
 the slide/swim moves in Fig.~4.
\smallskip

\begin{figure}[!h]
\includegraphics[scale=1.0]{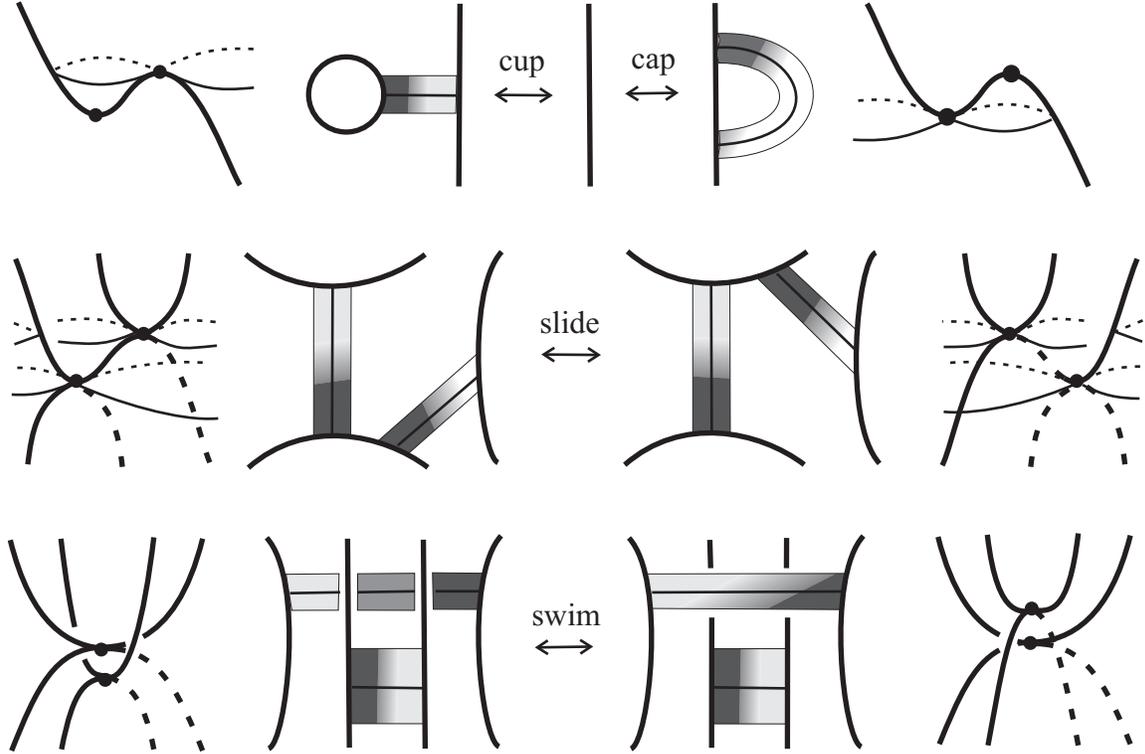}
\caption{Cup/cap moves and slide/swim moves of banded links}
\end{figure}

\noindent
(ii) If 2-links $S, S'$ are isotopic through generic 2-links,
 then the associated banded links $BL,BL'$ are related
 by the slide/swim moves in Fig.~4.
\smallskip

\noindent
(iii) If 2-links $S, S'$ are isotopic through generic 2-links 
 and one of $A_1^+A_1^+,A_1^+A_1^-,A_1^-A_1^-$-singularities,
 then $BL,BL'$ are related by the slide/swim moves in Fig.~4.
\smallskip

\noindent
(iv) If 2-links $S, S'$ are isotopic through generic 2-links 
 and exactly one $A_2$-singularity, then 
% the associated banded links 
 $BL,BL'$ are related
 by the cap/cup and slide/swim moves in Fig.~4.
\medskip

\noindent
\emph{Proof}.
(i)
The lowest critical point of a generic 2-link $S$ with respect to $\pr:S\to[-1,1]$ 
 at $t=t_1$ is a minimum, so the cross-section $S_{t_1+\ep}$ 
 is a trivial knot for some $\ep>0$.
The section $S_{t_1+\ep}$ is a prototype of a future banded link $BL$,
 which will be located in a fixed copy of $\R^3$.
The key idea in constructing $BL$ is to watch 
 the current cross-section $S_t=S\cap(\R^3\times\{t\})$
 simultaneously adding bands and trivial knots corresponding to 
 new saddles and minima, respectively.
The left column of Fig.~5 contains cross-sections $S_t$ for different values of $t$.
The right column shows successive stages of constructing $BL$ 
 whose final form is the top right.
\smallskip

\begin{figure}[!h]
\includegraphics[scale=1.0]{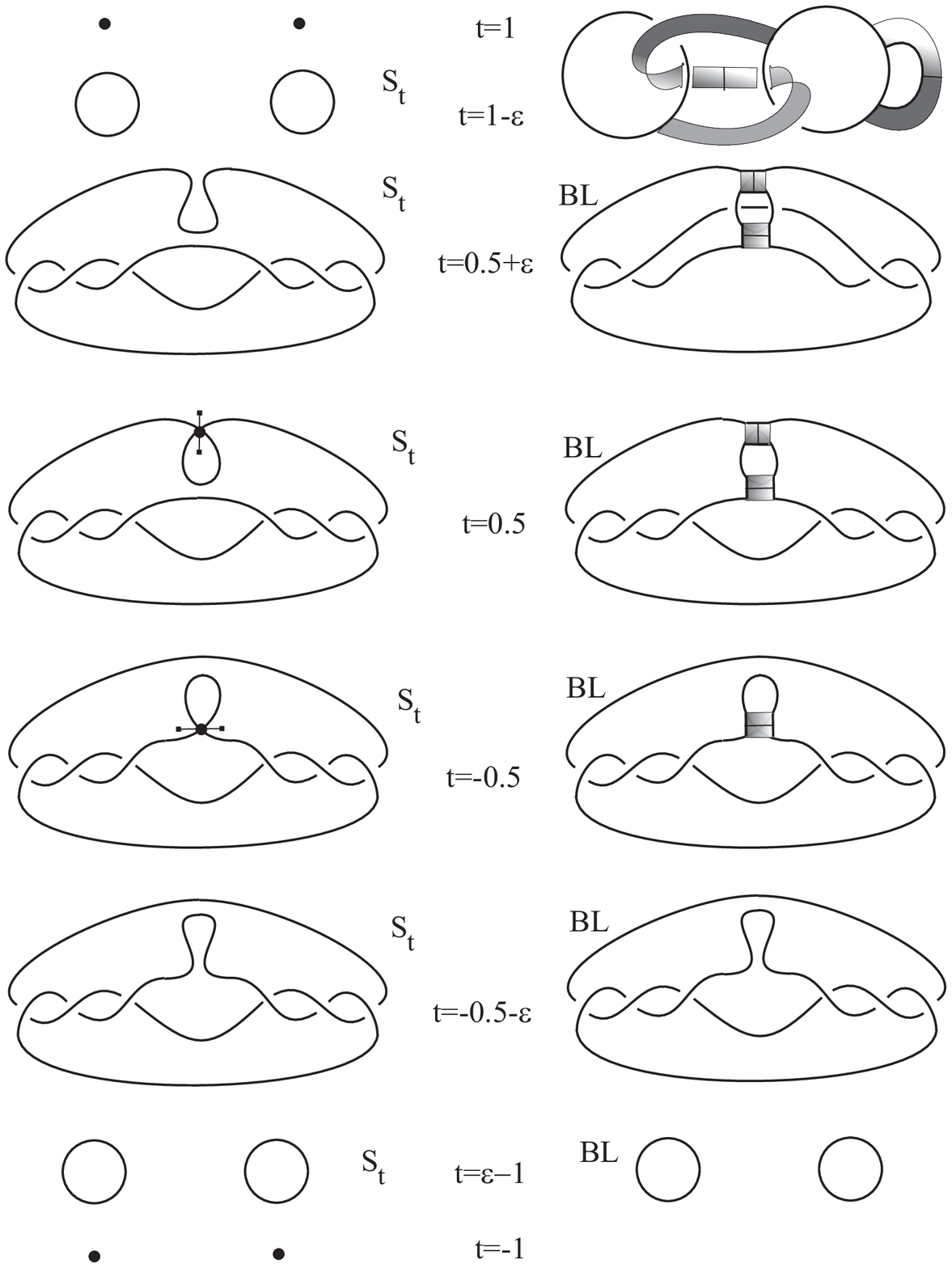}
\caption{Cross-sections and a banded link of the spun 2-knot of the trefoil}
\end{figure}

%\noindent
While $t$ is increasing, we isotopically deform 
 the current banded link $BL\subset\R^3$ 
 following $S_t=S\cap(\R^3\times\{t\})$, see Fig.~5.
The existing bands of $BL$ can be deformed 
 to avoid intersections with the rest of $BL$.
For each new minimum of $S$ in $\R\times\{t_j\}$, add a trivial knot 
 from $S_{t_j+\ep}$ to the current banded link $BL\subset\R^3$.
\smallskip

%\noindent
For each new saddle of $S$, attach a small band $B$ to $BL$.
The band $B$ has 2 opposite sides attached 
 to branches of the previous link $BL$.
While $t$ passes the critical value, the attached sides of $B$ are retracted 
 to a point and are replaced by the free sides of $B$.
The band $B$ can not meet the attached sides of other bands of $BL$
 since these sides are not included into the current cross-section of $S$.
So there are only 2 cases when the new link 
 with bands does not satisfy Definition~2.4.
\medskip

\noindent
(a) One (or two) of the attached sides of $B$ may meet 
 a free side of another band $B'$ of $BL$, see the upper picture of Fig.~6. 
Then slide $B$ along the free side of $B'$ in any of the two directions
 so that in the end the attached side of $B$ does not meet $B'$.
\medskip

\noindent
(b) The band $B$ intersects the interior of another band $B'$ of $BL$, 
 see the lower picture of Fig.~6.
Then $B$ swims through any of the attached sides of $B'$, so $B,B'$ fall apart.
The band $B$ can not swim through the free sides of $B'$ 
 as they belong to the current cross-section of $S$.
For each new 2-handle (a maximum), we keep the corresponding trivial knot 
 of $BL$, although it disappears from $S_t=S\cap(\R^3\times\{t\})$.
\medskip

After we have passed all critical values of $\pr:S\to\R$,
 the associated banded link $BL\subset\R^3$ has been constructed.
\medskip

\begin{figure}[!h]
\includegraphics[scale=1.0]{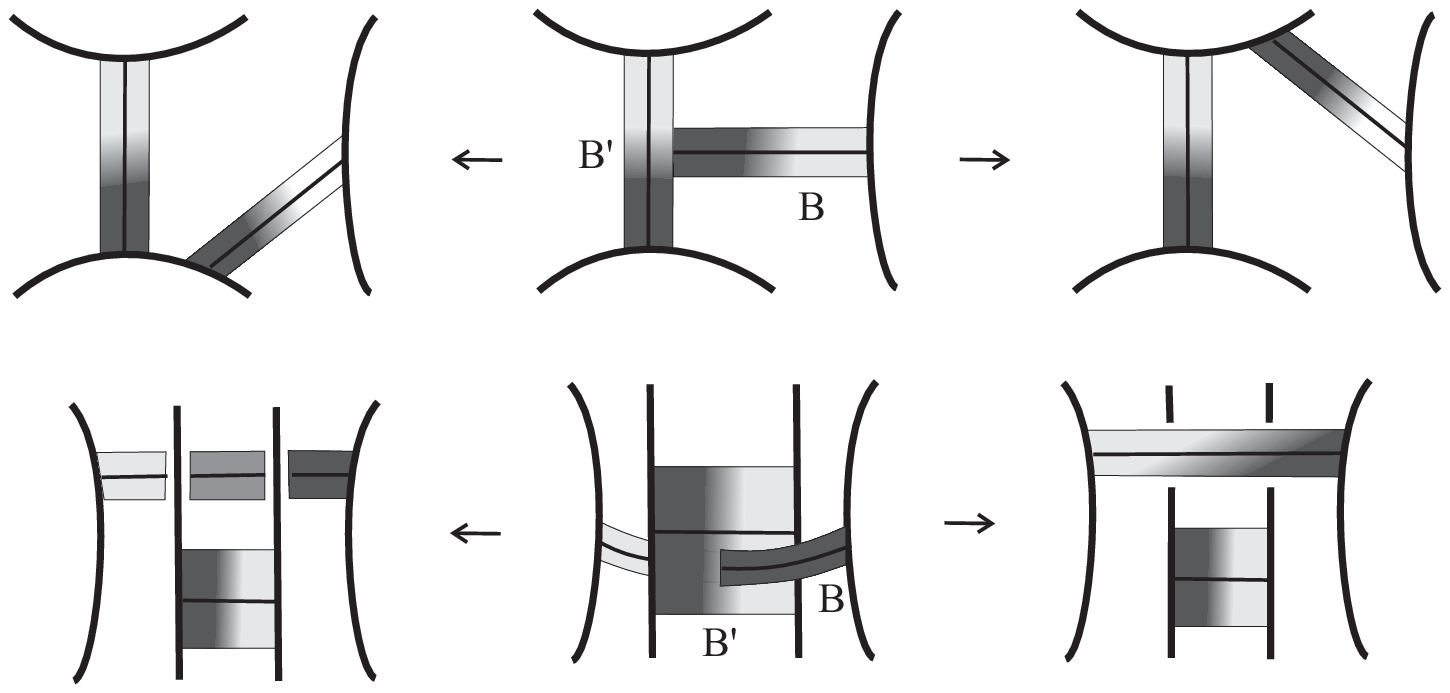}
\caption{A band slides or swims to remove an intersection}
\end{figure}
\medskip

\noindent
(ii) The construction above is not affected by an isotopy of $S$
 keeping the order of critical points of $\pr:S\to[-1,1]$.
Indeed all cross-sections $S_t$ are replaced by isotopic links,
 so the resulting banded link is isotopic to the original one
 provided that we remove intersections of bands in Fig.~6 in the same way.
\medskip

\noindent
(iii) 
The given isotopy of $S$ is a smooth path passing through
 one of $A_1^+ A_1^+, A_1^+ A_1^-, A_1^- A_1^-$-singularities in the space $\CS$ of 2-links.
For $A_1^+ A_1^+$ or $A_1^+ A_1^-$, 
 an extremum and another singularity swap their heights, so
 we add a new trivial knot (passing a minimum) or 
 keep an existing trivial knot (passing a maximum) that does not affect the other singularity.
For an $A_2$-singularity, two saddles of $S$ swap their heights, 
 so we add 2 bands to $BL$ in the reverse order.
Consider the critical moment when both saddles are 
 in the same section $\R^3\times\{t_j\}$.
If the associated bands do not intersect each other,
 then the new banded link is isotopic to the original one.
In (i) we listed the only cases (a), (b) when one band 
 may intersect another, which led to the moves in Fig.~6
 so the banded links are equivalent through
 the slide/swim moves.
\medskip

\noindent
(iv) 
If an isotopy of $S$ passes through an $A_2$-singularity, then 
 around this moment a non-degenerate saddle and extremum appear in a 2-link,
 see Claim~4.3(iv).
In the case of a minimum, one adds a trivial knot to the current banded link $BL$
 and a band attached to the trivial knot and to an existing branch of $BL$
 as shown in the cup move of Fig.~4.
In the case of a maximum, one adds a band attached by both sides to 
 a branch of the current banded link $BL$ as shown in the cap move of Fig.~4.
Recall that we keep the trivial knot when $t$ passes a maximum.
The leftmost and rightmost columns of Fig.~4 describe 
 projections of 2-links to $\R^3$ around singular moments.
The 4th axis of $\R^3\times\R$ projects to the vertical axis of $\R^3$.
\qed
\medskip

%\noindent
Conversely, any admissible banded link will give rise to 
 a 2-link in 4-space, see Lemma~2.8.
One can describe all moves of banded links 
 generating any isotopy of 2-links in 4-space.
Banded links were called \emph{knots with bands} in \cite{Swe}.

%------------------------------------------------------------------------------------------------

\subsection{Representing 2-links in 4-space by marked graphs in 3-space}%2.3
\noindent
\smallskip

\noindent
Theorem~1.4 is easier to prove representing 2-links by marked graphs, 
 which are singular links with bridges at singular points.
\medskip

\noindent
{\bf Definition 2.7.}
After deformation retracting each band of a banded link $BL$ to a point, 
 we get a \emph{singular} link \cite{KV}, a collection of closed curves
 with finitely many double transversal intersections, see Fig.~2, 4.
The core of each retracted band defines a bridge at the singular point,
 a straight arc in a small plane neighbourhood of each singular point.
We consider the resulting \emph{marked} graph $MG$ up to isotopy in $\R^3$
 keeping a neighbourhood of each singular point in a (moving) plane.
\medskip

\noindent
In the smooth approach, the zero section $S\cap(\R^3\times\{0\})$ 
 containing all saddles of $\pr:S\to[-1,1]$ is a marked graph
 whose bridges show how to resolve the singular points 
 for $t>0$ (along bridges) and $t<0$ (across bridges), see Fig.~2, 3.
An abstract marked graph $MG$, i.e. a singular link with bridges,
 can be converted into a banded link $BL$ replacing each bridge
 by a small rectangle whose core coincides with the bridge.
% and the opposite sides are attached to the branches of $MG$.
So there is a 1-1 correspondence between banded links and marked graphs.
Lemma~2.8 provides a unique function from the set of admissible banded links
 to the set of 2-links, which is the inverse of the function from Proposition~2.6.
\medskip

\noindent
{\bf Lemma 2.8.}
Any admissible banded link $BL\subset\R^3$ 
% (or, equivalently, any admissible marked graph $MG$) 
 gives rise to a 2-link $S\subset\R^4$ that 
 can be represented by $BL$ as in Proposition~2.6(i).
\smallskip

\noindent
\emph{Proof.}
Take the marked graph $MG\subset\R^3$ 
 associated to the given banded link.
Isotopically deform $MG$ in such a way that neighbourhoods 
 of all singular points of $MG$ are contained 
 in a single hyperplane of $\R^3\times\{0\}$.
\smallskip

Resolving the singular points along the bridges for $t>0$ 
 and across the bridges for $t<0$, extend the embedding 
 $MG\subset\R^3\times\{0\}$ to a surface $S'\subset\R^3\times[-\ep,\ep]$
 for some $\ep>0$, such that the boundary $\bd S'$ consists of 
 trivial links in $\R^3\times\{\pm\ep\}$. 
\smallskip

%\noindent
Since both sections $S'_{\pm\ep}=S'\cap(\R^3\times\{t=\pm\ep\})$ are unlinks,
 one can find isotopies $\ph_t^{\pm}:\R^3\to\R^3$, $t\in[\ep,1-\ep]$, such that
 each $\ph^{\pm}_{1-\ep}(S'_{\pm\ep})$ is 
 a collection of small disjoint circles in a plane.
The isotopies $\ph_t^{\pm}$ define the embedding of a 2-link $S$ without 
 small disks into $\R^3\times[\ep-1,1-\ep]$, one disk for each component of $\bd S$.
Attaching a disc to each boundary circle gives
 a closed surface $S\subset\R^3\times[-1,1]$.
\smallskip

The zero section $S\cap(\R^3\times\{0\})$ is the original marked graph $MG$.
A small isotopy deformation makes $S$ generic.
The construction of Proposition~2.6(i) gives a banded link 
 equivalent to $MG$ as all bands may be chosen small and non-intersecting.
\qed
\medskip

%===================================================

\section{Three-page embeddings of marked graphs}%3

%--------------------------------------------------------------------------

\subsection{Any marked graph can be embedded into the 3-page book}%3.1
\noindent
\smallskip

\noindent
Recall that the 3-page book is $\tb=\R\times T$, where
 $T$ is the triod consisting of 3 edges $E_0,E_1,E_2$
 joining the vertex $O$ to the other 3 vertices.
The line $\al=\R\times O$ is said to be the \emph{binding} axis,
 $P_i=\R\times E_i$ are called the \emph{pages}, $i=0,1,2$.
\medskip

\begin{figure}[!h]
\includegraphics[scale=1.0]{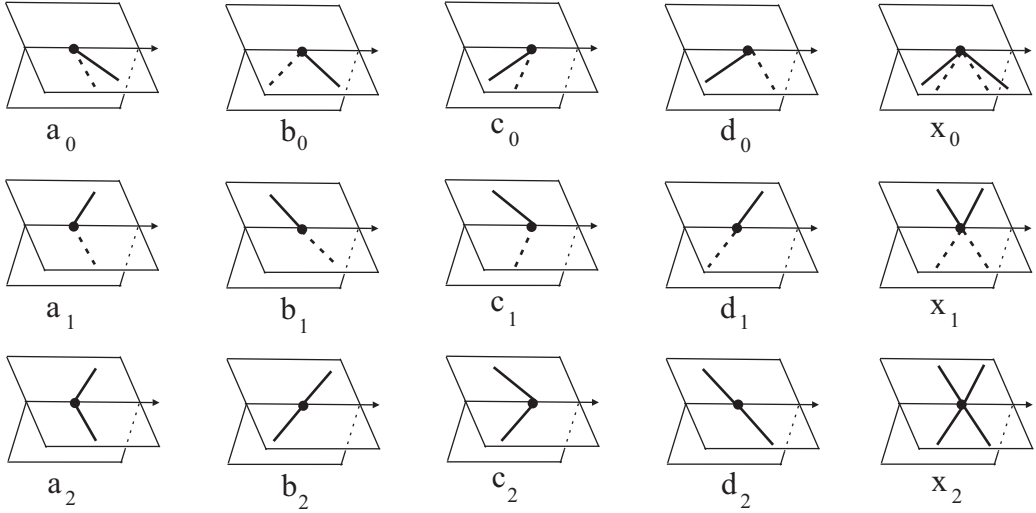}
\caption{The encoding letters for 3-page embeddings of marked graphs}
\end{figure}

\noindent
{\bf Definition 3.1.}
An embedding of a marked graph $G$ into the 3-page book $\tb$
 is called a \emph{3-page} embedding, if conditions (i)-(v) hold:
\smallskip

\noindent
(i) the intersection $G\cap\al$ of $G$ and the binding axis $\al$ is a finite set of points;
\smallskip

\noindent
(ii) the arcs at every point of $G\cap\al$ lie in 2 pages $P_i,P_j$, 
 $i\neq j$, see Fig.~7;
\smallskip

\noindent
(iii) all singular points of $G$ lie in $\al$,
 a neighbourhood of each singular point lies\\
\hspace*{7mm} 
 in a broken plane consisting of two pages 
 and looks locally like a cross $\times$;
\smallskip

\noindent
(iv) the bridge at each singular point lies in the binding axis $\al$;
\medskip

\noindent
(v) every connected component of $G\cap P_i$ is projected
 monotonically to $\al$.
\bigskip

\noindent
The arcs in the page $P_2$ are dashed in Fig.~7, 8.
All classical and singular links can be embedded into $\tb$ 
 in the sense of Definition~3.1, see Fig.~8.
\medskip

The pictures in each vertical column of Fig.~7 are obtained
 from each other by rotation around $\al$.
The rotation corresponds to the shift $i\mapsto i+1$ of indices,
 $i\in\Z_3=\{0,1,2\}$.
A 3-page embedding can be encoded by a word in the alphabet
 of 15 letters describing the local behaviour of $G$
 near the intersection points $G\cap\al$, see Fig.~7.
The 3-page embedding in Fig.~8 is encoded by 
 $w_G=a_0 a_1 (b_2b_0b_1)^2 d_0 a_1 
(x_1 b_1)^2 c_1 d_1 b_0 (d_1 d_0 d_2)^2 c_1 c_0$.
So a 3-page embedding of the marked graph $G_S$ of a 2-link $S$
 is a 1-dimensional representation of $S\subset\R^4$.
\medskip

\begin{figure}[!h]
\includegraphics[scale=1.0]{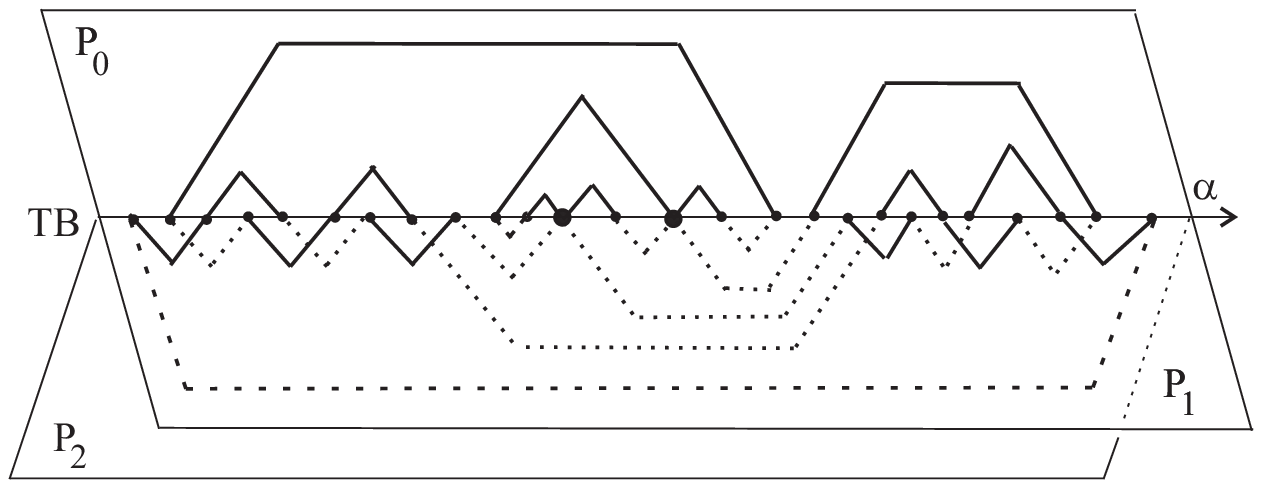}
\caption{A 3-page embedding of the marked graph from Fig.~9}
\end{figure}
%\medskip

We give a proof of the embedding result from \cite{KV},
 because this construction plays an important role
 in further considerations.
\medskip

\noindent
{\bf Proposition~3.2.} \cite{KV}
Any marked graph $G\subset\R^3$ is isotopic to 
 a 3-page embedding $G\subset\tb$ in the sense of Definition~3.1.
\medskip

\noindent
\emph{Proof}.
Consider a plane diagram $D$ of $G\subset\R^3$ in general position 
 with finitely many double crossings.
At each crossing in the diagram $D$ mark a small overcrossing arc.
Recall that, at each singular point of $G$, there is a marked bridge 
 transversally intersecting both branches of $G$ passing through the singular point.
\smallskip

%\noindent
In the plane containing the diagram $D$, draw a continuous path $\al$ such that
\smallskip

\noindent
(1) the path $\al$ passes through each marked arc and bridge exactly once;
\smallskip

\noindent
(2) $\al$ transversally intersects the rest of $D$, 
 the endpoints of $\al$ are away from $D$.
\medskip

\begin{figure}[!h]
\includegraphics[scale=1.0]{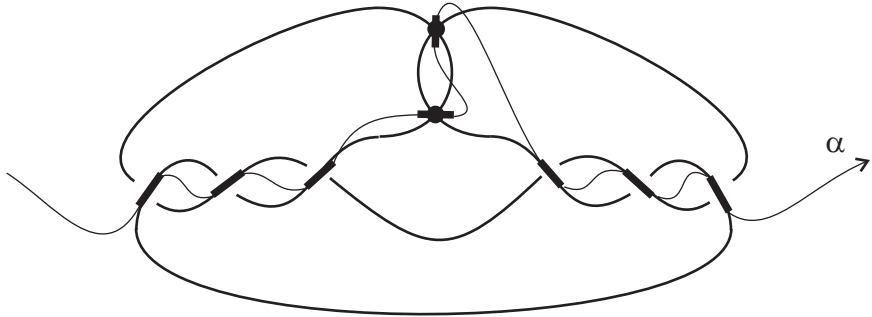}
\caption{How to construct a 3-page embedding of a marked graph}
\end{figure}

Isotopically deform the plane containing $D$ in such a way that $\al$ becomes
 a straight line containing all marked arcs and bridges of $D$.
Denote the upper half-plane and lower half-plane of $\R^2-\al$ 
 by $P_0$ and $P_2$, respectively.
Notice that a neighbourhood of each singular point 
 looks like a cross $\times$ with a centre in the axis $\al$, see Fig.~9.
\smallskip

Attach the third half-plane $P_1$ to $\al$ and 
 push all marked arcs into $P_1$, see Fig.~8.
If both (say) upper arcs at some singular point $v\in G$ 
 go to points on one side of the point $v\in\al$,
 then make an additional couple of crossings in the intersection $\al\cap D$ like
 in Reidemeister move II, see Fig.~13.
For instance, in the embedding $a_2 b_2 x_2$ both upper arcs go to the right,
 see the lower right picture of Fig.~15 \cite{KV}.
Then the intersection $G\cap P_i$
 is a finite collection of disjoint arcs, which can be made monotonic
 with respect to the projection $\tb\to\al$, $i=0,1,2$.
\qed
\smallskip

%--------------------------------------------------------------------------

\subsection{Any isotopy of links can be realised in the hexabasic book}%3.2
\noindent
\smallskip

The following lemma is a key stone of the 3-page approach to knot theory.
\medskip

\noindent
{\bf Lemma 3.3.} \cite{Kur}
Any isotopy of 3-page embeddings of classical links is decomposed
 into finitely moves in Fig.~10 and theirs images under
 $i\mapsto i+1$, $i\in\Z_3$. 

\begin{figure}[!h]
\includegraphics[scale=1.0]{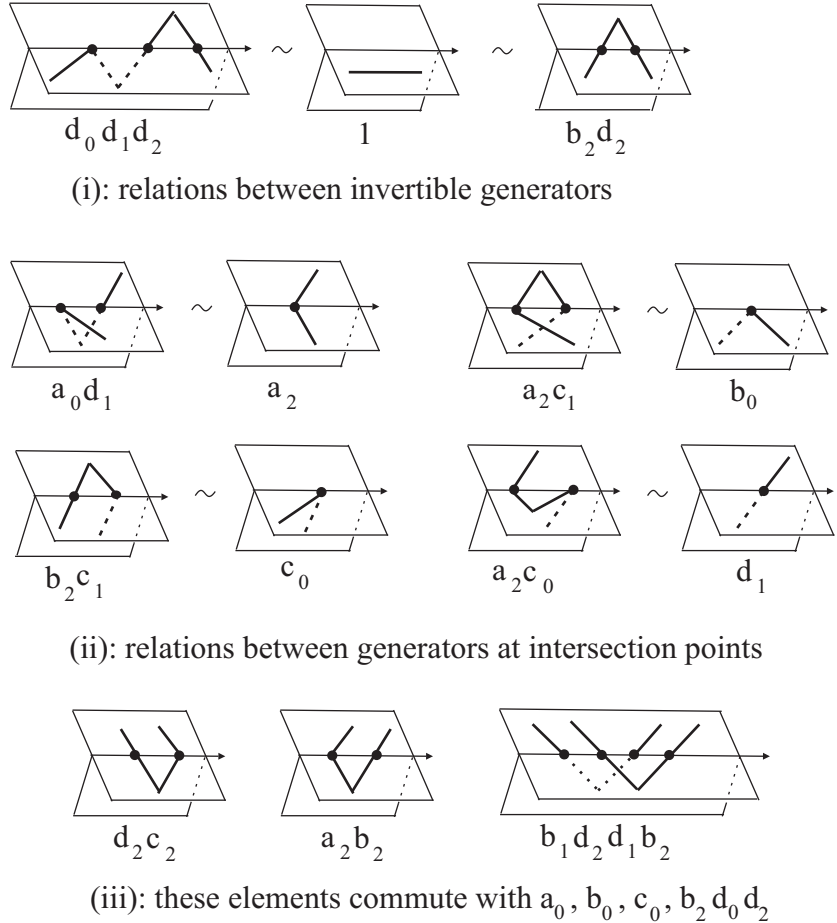}
\caption{Finitely many moves generating any isotopy of classical links}
\end{figure}
\medskip

%\noindent
The algebraic form of the moves in Fig.~10 is below,
  $i\in\Z_3=\{0,1,2\}$, see \cite{Kur}:
\smallskip

\noindent
(1) $d_0d_1d_2=1$, $b_id_i=1=d_ib_i$;
\smallskip

\noindent
(2) $a_i=a_{i+1}d_{i-1}$, $b_i=a_{i-1}c_{i+1}$, 
     $c_i=b_{i-1}c_{i+1}$, $d_i=a_{i+1}c_{i-1}$;
\smallskip

\noindent
(3) $uv=vu$, where $u\in\{a_ib_i,d_ic_i,b_{i-1}d_id_{i-1}b_i\}$,
     $v\in\{a_{i+1},b_{i+1},c_{i+1},b_id_{i+1}d_i\}$.
\medskip

Lemma~3.4 is the crucial step in Theorem~1.3.
\medskip

\noindent
{\bf Lemma 3.4.}
The moves in Fig.~10 are realised in the hexabasic book $\hb$.
\smallskip

\noindent
\emph{Proof.} All the moves in Fig.~10, apart from the commutativity
 of $a_i,b_i,c_i,b_{i-1}d_id_{i-1}$ with $b_{i+1}d_{i-1}d_{i+1}b_{i-1}$,
 can be realised in the 3-page book $\tb$.
For instance, the relation $b_2d_2=1$ is realised
 by compressing the slice between the 2 intersection points
 and removing the resulting point from $\al$.
The other relations are realised in $\hb$, see a geometric realisation
 of $(b_1d_2d_1b_2)a_0=a_0(b_1d_2d_1b_2)$ in Fig.~11.
\qed
\medskip

\begin{figure}[!h]
\includegraphics[scale=1.0]{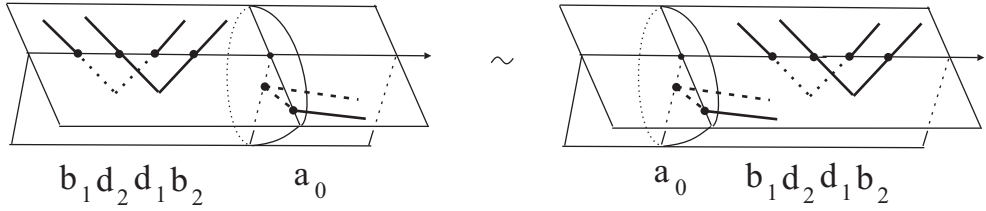}
\caption{Realising a commutative relation in the hexabasic book $\hb$}
\end{figure}

%--------------------------------------------------------------------------

\subsection{Any 2-link is isotopic to a surface in the universal polyhedron}%3.3
\noindent
\smallskip

Here we finish the proof of Theorem~1.3.
\medskip

\noindent
\emph{Proof of Theorem~1.3}.
By Lemma~2.3 any 2-link in 4-space is isotopic to 
 a surface $S\subset\R^3\times[-1,1]$ having
 all maxima, minima and saddles in
 different sections $\R^3\times\{t_j\}$ 
 for some $-1<t_1<\dots<t_n<1$.
In Step~1 we embed each cross-section $S_{t_j}$ into the 3-page book.
In Step~2 we extend this embedding to a regular neighbourhood of $S_{t_j}$. 
In Step~3 we embed the complement of the neighbourhoods into $\hb\times[-1,1]$.
\smallskip

\noindent
\emph{Step 1}.
Choose $\ep>0$ such that the closed $\ep$-neighbourhoods 
 $N_{\ep}(S_{t_j})$ of $S_{t_j}$ in $S$ are disjoint and each of them
 contains exactly one critical point of $\pr:S\to[-1,1]$, $j=1,\dots,n$.
Then the boundaries $\bd N_{\ep}(S_{t_j})$ are classical links.
By Proposition~3.2 there is an isotopy $f^u_j:\R^3\times\{t_j\}\to\R^3\times\{t_j\}$, 
 $u\in[0,1]$, moving $S_{t_j}$ into $\tb$, i.e. $f^0_j=\id_{\R^3}$,
 $f^1_j(S_{t_j})\subset\tb\times\{t_j\}$ is a 3-page embedding.
Take smooth functions $g_j:[t_j-\ep,t_j+\ep]\to[0,1]$ 
 such that $g_j(t_j)=1$ and $g_j(t_j\pm\ep)=0$.
Extend $f_u^j$ to
 $$F^u_j:\R^3\times[t_j-\ep,t_j+\ep]\to\R^3\times[t_j-\ep,t_j+\ep],
     \quad u\in[0,1],$$ 
 $$F^u_j(x,t)=(f^{ug_j(t)}_j(x),t), \mbox{ where }x\in\R^3, t\in[t_j-\ep,t_j+\ep].$$
Then $F^u_j=f^u_j$ for $t=t_j$ and $F^u_j=\id$ for $t=t_j\pm\ep$.
Hence $\bd N_{\ep}(S_{t_j})$ are pointwise fixed and we may 
 combine $F^u_j$ together to form a smooth isotopy 
 $F^u:\R^3\times[-1,1]\to\R^3\times[-1,1]$ moving
 each $S_{t_j}$ into $\tb\times\{t_j\}$.
Denote the resulting surface by $S'$.
\medskip

\noindent
\emph{Step 2}.
If a singular cross-section $S'_{t_j}$ has a double intersection, then
 both positive and negative resolutions of $S'_{t_j}$ can be embedded into $\tb$.
Indeed the positive and negative resolutions of the singular point $x_i$
 are encoded by $1$ and $c_ia_i$, respectively, see Fig.~12.
Given an encoding word $w_j$ of $S'_{t_j}\subset\tb$,
 the positive resolution of $S'_{t_j}$ is encoded by $w_j$ after removing 
 the letter $x_i$ representing the double point of $S'_{t_j}$.
\smallskip

\begin{figure}[!h]
\includegraphics[scale=1.0]{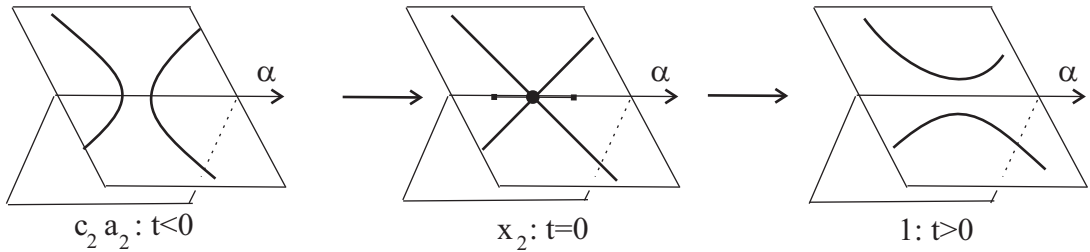}
\caption{Resolving a singular point in the 3-page book $\tb$}
\end{figure}

The argument below with the sign $\pm$ covers 2 cases 
 when either $+$ or $-$ is taken in all formulae.
If $S'_{t_j}$ contains a maximum or minimum, 
 $S'_{t_j\pm\ep/2}$ can be embedded into $\tb$.
So there are isotopies 
 $h^u_{\pm j}:\R^3\times\{t_j\pm\ep/2\}\to\R^3\times\{t_j\pm\ep/2\}$, $u\in[0,1]$,
 moving each $S'_{t_j\pm\ep/2}$ into $\tb\times\{t_j\pm\ep/2\}$.
Take smooth functions $\ti g_{j}:[t_j-\ep,t_j+\ep]\to[0,1]$ such that 
 $\ti g_j(t_j)=0=\ti g_j(t_j\pm\ep)$ and $\ti g_j(t_j\pm\ep/2)=1$.
Extend $h^u_{\pm j}$ to
$$H^u_{j}:\R^3\times[t_j-\ep,t_j+\ep]\to\R^3\times[t_j-\ep,t_j+\ep],
    \quad u\in[0,1],$$ 
$$H^u_{j}(x,t)=(h^{u\ti g_j(t)}_{\pm j}(x),t) \mbox{ for } x\in\R^3,\;
 t \mbox{ between } t_j \mbox{ and } t_j\pm\ep.$$
Then $H^u_j=h^u_{\pm j}$ for $t=t_j\pm\ep/2$ and 
 $H^u_j=\id$ for $t=t_j$, $t=t_j\pm\ep$.
Hence $S'_{t_j}$ and $\bd N_{\ep}(S'_{t_j})$ are pointwise fixed and 
 we may combine $H^u_{j}$ together to form a smooth isotopy 
 $H^u:\R^3\times[-1,1]\to\R^3\times[-1,1]$ moving
 each $N_{\ep/2}(S'_{t_j})$ into $\tb\times[t_j-\ep/2,t_j+\ep/2]$.
Denote the resulting surface by $S''$.
\medskip

\noindent
\emph{Step 3}.
The cross-sections $S''_{t_j+\ep/2}$ and $S''_{t_{j+1}-\ep/2}$
 are isotopic classical links, $j=1,\dots,n-1$.
By Lemma~3.3 and Proposition~3.4 any isotopy of classical links
 can be realised in $\hb$.
Then the layers $S''\cap(\R^3\times[t_j+\ep/2,t_{j+1}-\ep/2])$
 can be replaced by an isotopy of links in $\hb\times[t_j+\ep/2,t_{j+1}-\ep/2]$.
It remains to extend the embedding to the neighbourhoods
 of the lowest minimum and highest maximum of $S''$
 shrinking their boundaries in $\hb$.
So the final surface is embedded into $\hb\times[-1,1]$.
\qed
\medskip

%===================================================

\section{The universal semigroup of 2-dimensional links}%4

%------------------------------------------------------------------------------------------------

\subsection{Local moves of marked graphs generate any isotopy of 2-links}%4.1
\noindent
\smallskip
 
Here we derive a complete set of moves of banded links and marked graphs,
 that generate any isotopy of 2-links in 4-space.
Marked graphs can be represented by plane diagrams with 
 small straight arcs denoting bridges over singular points, see Fig.~2, 9.
In particular, the cyclic order of edges at each singular point is invariant.
\medskip

\noindent
{\bf Lemma~4.1} \cite{Kau}
Marked graphs are isotopic in $\R^3$ if and only if
 their plane diagrams can be obtained from each other
 by finitely many Reidemeister moves in Fig.~13,
 where all symmetric images of the moves should be considered.
\medskip

\begin{figure}[!h]
\includegraphics[scale=1.0]{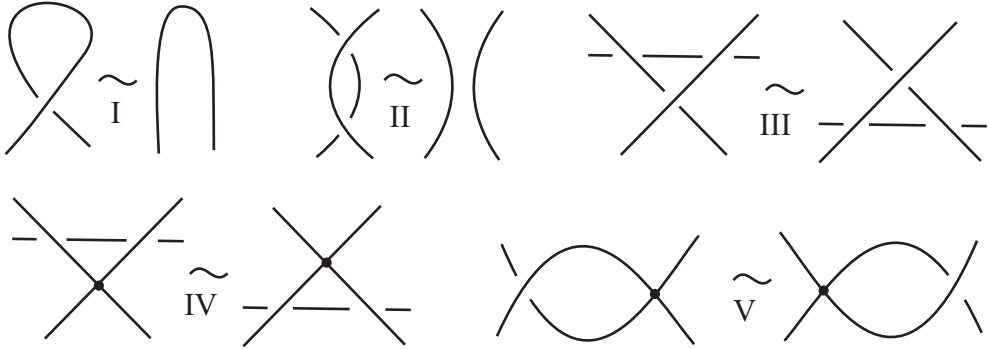}
\caption{Reidemeister moves for rigid isotopy of marked graphs}
\end{figure}

The moves in Fig.~13 are \emph{local} in the sense,
 that a small disk in the left part of each move is replaced by 
 another small disk in the right part of the move, while
 the rest of the diagram remains unchanged.
The singular points in moves IV and V of Fig.~13 
 can be equipped with arbitrary corresponding bridges.
The proof is a direct application of the transversality theorem of Thom
 similarly to a proof of the Reidemeister theorem 
 for plane diagrams of classical links, see \cite[section~2]{FK}.
\medskip

\noindent
{\bf Proposition~4.2.} 
%(Yoshikawa \cite{Yos})
Marked graphs represent isotopic 2-links in 4-space
 if and only if they can be obtained from each other 
 by finitely many moves in Fig.~14.
\medskip

\begin{figure}[!h]
\includegraphics[scale=1.0]{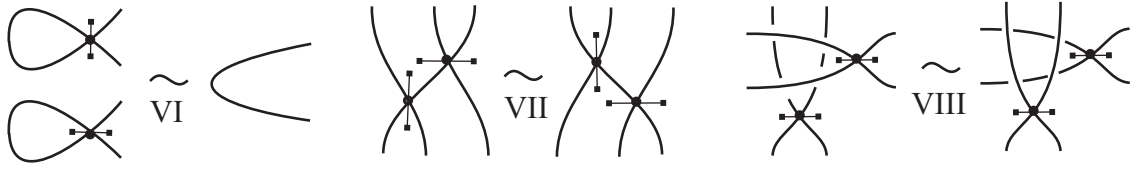}
\caption{Moves of marked graphs generating isotopy of 2-links}
\end{figure}

Symmetric images of the moves in Fig.~14 are skipped as 
 they can be reduced to the standard moves using an isotopy in $\R^3$.
Proposition~4.2 was conjectured by K.~Yoshikawa in \cite{Yos}.
F.~Swenton \cite{Swe} claimed a proof of Proposition~4.2 
 using banded links and the equivalent moves in Fig.~4.
M.~Saito wrote in his review for the MathSciNet:
`It is claimed that this set of moves is equivalent to Yoshikawa's moves.
 It might be beneficial of some more detailed accounts, for example,
 those for the above claim, are discussed further and 
 presented elsewhere in the literature'.
The authors were asked by S.~Carter to fill in these details,
 so we give a more detailed proof of Proposition~4.2 for banded links.
Recall that the singular subspace $\Si$ of the space $\CS$ 
 of 2-links was introduced in Definition~2.5.
The following result will be formally deduced in Appendix using
 the transversality theorem of Thom.
\medskip

\noindent
{\bf Claim~4.3.}
(i) The closure of the subspace $\Si$ has codimension~1 in the space $\CS$.
\smallskip

\noindent
(ii)
The complement of the closure $\overline{\Si}$ in $\CS$ consists of generic 2-links.
\smallskip

\noindent
(iii)
Any isotopy of 2-links can be deformed in such a way that 
 all intermediate 2-links are generic apart from 
 finitely many singularities of Definition~2.5.
\smallskip

\noindent
(iv)
If an isotopy passes through an $A_2$-singularity, then 
 a non-degenerate saddle and extremum collide and disappear
 as shown in the top picture of Fig.~18.
\medskip

\noindent
Claims~4.3(i,ii) say that any point of $\CS$ can be removed from $\Si$
 by a small perturbation, i.e. a 2-link can be made generic, which implies Claim~2.3(ii).
Claim~4.3(iii) says that the singularities of Definition~2.5 are the only singularities
 occuring in any isotopy of 2-links in general position.
\medskip

\noindent
\emph{Proof} of Proposition~4.2.
By Claim~4.3(iii) any isotopy of 2-links can be deformed into a smooth path 
 transversal to the subspace $\Si\subset\CS$.
When the path passes through one of the singularities, the associated 
 banded link changes according to Proposition 2.6(iii),(iv),
 which led to the moves in Fig.~4 as required.
\qed
\medskip

%----------------------------------------------------

\subsection{A 1-dimensional encoding 2-links up to isotopy in 4-space}%4.2
\noindent
\smallskip

Here we reduce the isotopy classification of 2-links in 4-space 
 to a word problem in the finitely presented semigroup $\SL$,
 the universal semigroup of 2-links.
Recall that moves (1)-(8) on 3-page embeddings were defined in subsection~1.3.
Theorem~1.4 follows from the following generalisation of 
 Lemma~3.3 to singular links.
\medskip

\noindent
{\bf Proposition~4.4.} \cite{KV}
Consider the semigroup $\SK$ generated by $a_i,b_i,c_i,d_i,x_i$, 
 $i\in\Z_3$, subject to relations (1)-(5) from subsection~1.3.
Then any singular link $G\subset\R^3$ is encoded by an element
 $w_G\in\SK$ in such a way that singular links $G,G'$ are isotopic
 if and only if their encoding elements $w_G$ and $w_{G'}$
 are equal in $\SK$.
An element $w\in\SK$ encodes a singular link 
 if and only if $w$ is central in $\SK$.
\medskip

\noindent
\emph{Proof of Theorem~1.4.}
Any 2-link can be represented by its marked graph $G$ whose 3-page embedding
 is encoded by a word in the letters $a_i,b_i,c_i,d_i,x_i$, $i\in\Z_3$,
 as described before Proposition~3.2.
All encoding elements form the centre of $\SL$ as the same result
 holds for the universal semigroup $\SK$ of singular links, i.e. relations 
 (1)-(5) imply that any encoding element commutes with the generators.
\smallskip

The remaining part of Theorem~1.4 states that two 3-page embeddings 
 of marked graphs represent isotopic 2-links in 4-space if and only if
 they can be related by algebraic moves (1)--(8) in subsection~1.3.
By Lemmas~3.3-3.4 and Proposition 4.2 it suffices to realise 
 moves VI, VII, VIII in Fig.~14 by 3-page embeddings.
\smallskip

In moves VI, VII, VIII a small disk in the left part 
 is replaced by another small disk in the right part.
Similarly to the construction of a 3-page embedding, choose a path $\al$
 passing through overcrossing arcs and bridges at singular points, see Fig.~15, 16, 17.
Deform the diagrams in such a way that $\al$ becomes a straight line and 
 push all overcrossing arcs into the half-plane $P_1$, all bridges remain in $\al$.

\begin{figure}[!h]
\includegraphics[scale=1.0]{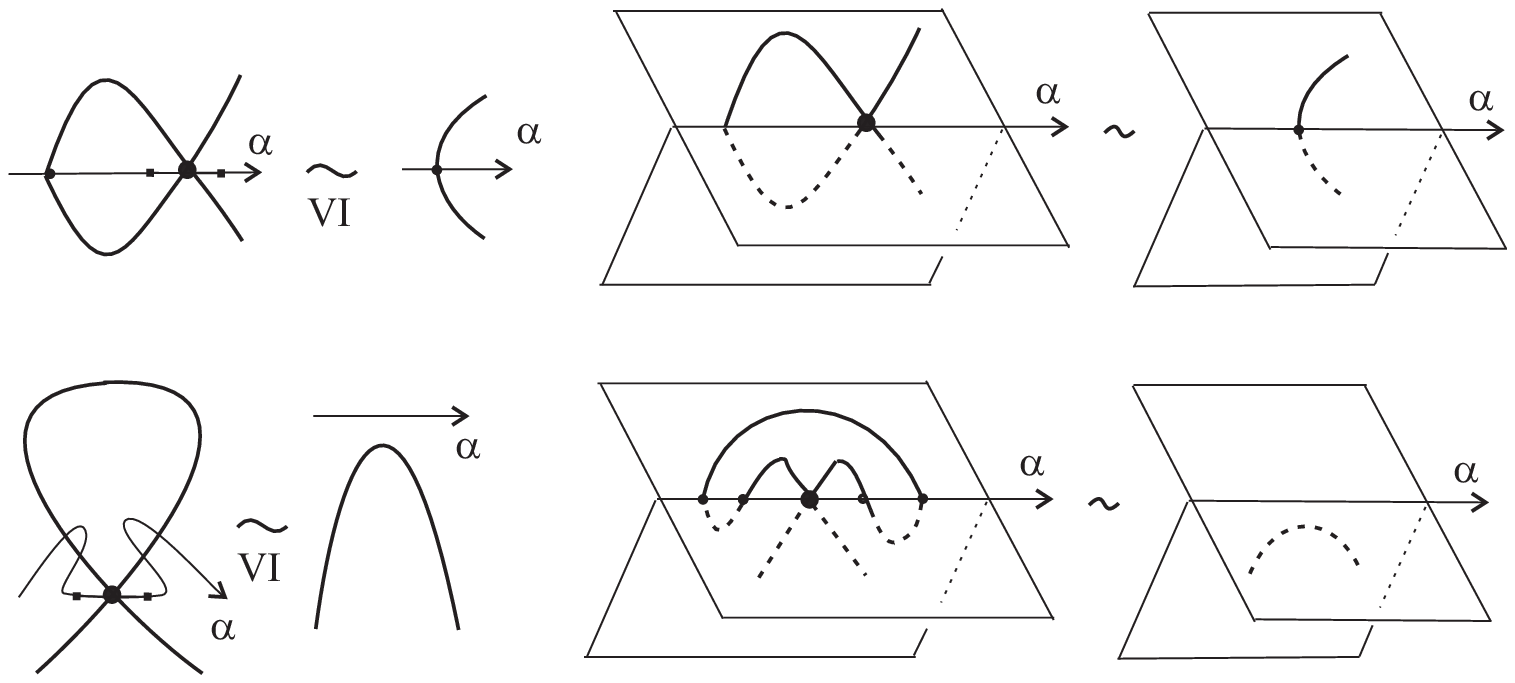}
\caption{Realising moves VI of Fig.~14 in terms of 3-page embeddings}
\end{figure}

In Fig.~15 moves VI are encoded by $a_1 x_1=a_1$ and 
 $a_1 b_1 x_1 d_1 c_1=1$ equivalent to  (6) for $i=1$.
We made additional intersections of $\al$ with the diagram to decompose the resulting
 embedding into local 3-page embeddings from Fig.~7.

\begin{figure}[!h]
\includegraphics[scale=1.0]{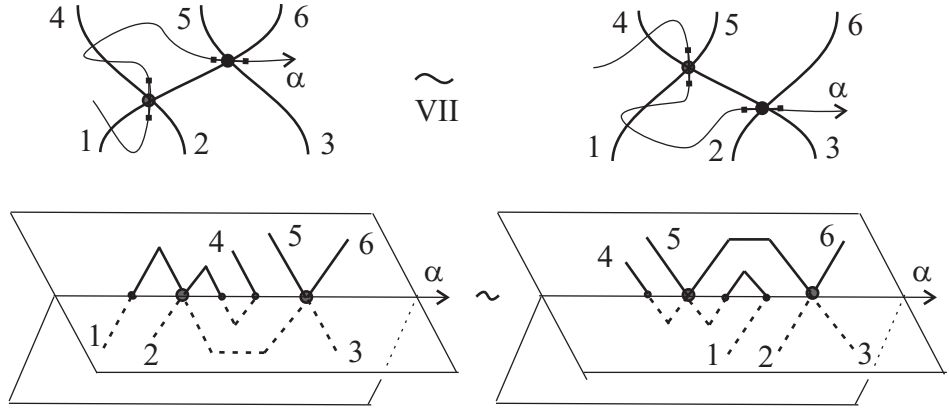}
\caption{Realising move VII of Fig.~14 in terms of 3-page embeddings}
\end{figure}

In Fig.~16 move VII is encoded by 
 $d_1 x_1 b_1 c_1 x_1=b_1 x_1 d_1 c_1 x_1$, 
 which is (7) for $i=1$.
Numbers 1, 2, 3, 4, 5, 6 denote arcs going out of the small disk replaced by move VII,
 e.g. the path $\al$ starts between arcs 1, 4 and ends between arcs 3, 6. 

\begin{figure}[!h]
\includegraphics[scale=1.0]{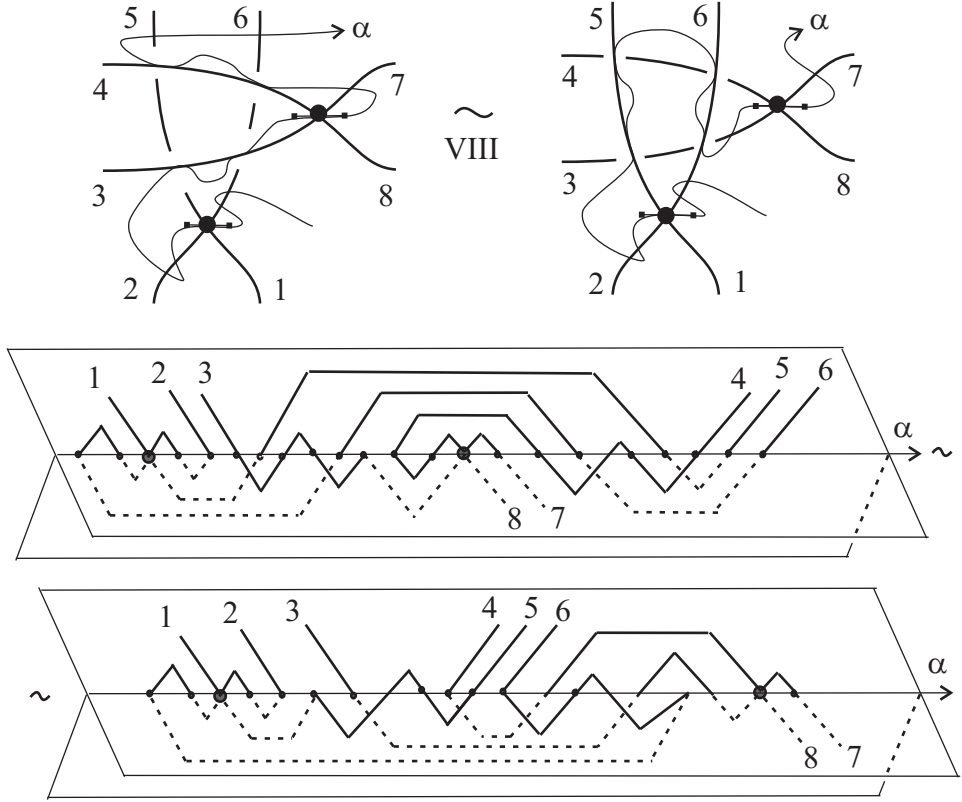}
\caption{Realising move VIII of Fig.~14 in terms of 3-page embeddings}
\end{figure}

In Fig.~17 move VIII is encoded by
 $$(a_1 b_1 x_1 b_1 c_1) d_2 d_1 (b_2 d_2) 
 d_1 d_0 a_2 b_2 x_1 b_1 d_2 b_1 (b_2 d_2) b_1 b_2 d_1^2=$$
 $$(a_1 b_1 x_1 b_1 c_1) b_0 b_1 (b_2 d_2) a_1 b_2 a_2 d_1 (b_2 d_2) d_1 c_0 b_1 x_1 b_1,$$
 which is equivalent to (8) for $i=1$ after removing $b_2 d_2=1$ by relation~(1).
The relations for other $i\in\Z_3$ were added 
 to make the presentation symmetric.
\qed 

%==================================================================

\section*{Appendix: the multi-jet transversality theorem of Thom}%A
\smallskip

Here we follow \cite[sections I.2, I.8]{AVG}.
Let $\xi,\eta:M\to N$ be smooth maps between
 finite dimensional manifolds
 with Riemannian metrics $\rho_M,\rho_N$, respectively.
\medskip

\noindent
{\bf Definition~A.1.}
The maps $\xi$ and $\eta$ have the tangency of \emph{order} $k$
 at a point $z\in M$ if $k$ is the maximal integer such that
 $\rho_N(\xi(w),\eta(w))/(\rho_M(z,w))^k\to 0$ as $w\in M$ tends to $z$, e.g.
 the curve $\xi(w)=w^{k+1}$ has the tangency of order $k$ with $\eta(w)=0$. 
\smallskip

%\noindent
The $l$-tuple $k$-\emph{jet} of the map $\xi$ at
 $(z_1,\dots,z_l)\in M^l$ is the equivalence class
 of smooth maps $\eta:M\to N$ up to tangency of order $k$
 at the points $z_1,\dots,z_l\in M$, e.g.
 the 1-tuple 1-jet $j^1_{[1]}\xi(z)$ of a map $\xi:\R\to\R$
 is determined by $z,\xi(z),\dot \xi(z)$.
\medskip

Denote by $J^k_{[l]}(M,N)$ the space of all $l$-tuple $k$-jets of
 smooth maps $\xi:M\to N$ for all $(z_1,\dots,z_l)\in M^l$.
Let $(x_1,\dots,x_m)$ and $(y_1,\dots,y_n)$ be
 local coordinates in $M$ and $N$, respectively.
If $\xi$ is defined locally by
 $y_j=\xi_j(x_1,\dots,x_m)$, $j=1,\dots,n$, then
 the $l$-tuple $k$-\emph{jet} of $\xi$ at $(z_1,\dots,z_l)$
 is determined by $l$ arrays of the data below
$$\{x_1,\dots,x_m\};\quad
  \{y_1,\dots,y_n\};\quad
  \left\{ \dfrac{ \bd\xi_j }{ \bd x_i } \right\};\quad
  \; \ldots \;
  \left\{ \dfrac{ \bd^{k}\xi_j }{ \bd x_{i_1}\dots x_{i_s} } \right\},
  i_1+\dots +i_s=k.$$

The quantities above define local coordinates
 in $J^k_{[l]}(M,N)$.
The $l$-tuple $k$-\emph{jet} $j^k_{[l]}\xi$ of
 a smooth map $\xi:M\to N$ can be considered as
 the map $j^k_{[l]}\xi:M^l\to J^k_{[l]}(M,N)$,
 namely $(z_1,\dots,z_l)$ goes to
 the $l$-tuple $k$-jet of $\xi$ at $(z_1,\dots,z_l)$.
\medskip

The manifold $J^k_{[l]}(M,N)$ is finite dimensional, e.g.
 $J^0_{[l]}(M,N)=(M\times N)^l$,
$$\dim J^1_{[l]}(M,N)=(m+n+mn)l, \;
    \dim J^2_{[l]}(M,N)=(m+n+mn+\frac{m(m+1)}{2}n)l.$$
\smallskip

\noindent
{\bf Definition~A.2.}
Take an open set $W\subset J^k_{[l]}(M,N)$.
The set of smooth maps $f:M\to N$ with $l$-tuple $k$-jets from $W$
 is \emph{open}.
These sets for all open $W\subset J^k_{[l]}(M,N)$ over all $k,l$
 form a basis of the \emph{Whitney} topology in $C^{\infty}(M,N)$.
The space $\CS$ of all 2-links $S\subset\R^4$ inherits
 the \emph{Whitney} topology from $C^{\infty}(S,\R^4)$.
So two maps are close in the Whitney topology
 if they are close with all derivatives.
\medskip

\noindent
{\bf Definition A.3.}
Let $M$ be a finite dimensional smooth manifold.
A subspace $\La\subset M$ is called
 \emph{a stratified space} if $\La$ is the union
 of disjoint smooth submanifolds $\La^i$ (\emph{strata})
 such that the boundary of each stratum
 is a finite union of strata of less dimensions.
Let $N$ be a finite dimensional manifold.
A smooth map $\xi:M\to N$ is \emph{transversal} to
 a smooth submanifold $U\subset N$ if
 the spaces $\xi_*(T_zM)$ and $T_{\xi(z)}U$ generate $T_{\xi(z)}N$
 for each $z\in M$.
A smooth map is $\eta:M\to V$ \emph{transversal} to
 a stratified space $\La\subset V$
 if the the map $\eta$ is transversal to each stratum of $\La$.
\medskip

\noindent
Briefly Theorem~A.4 says that any map
 can be approximated by `a nice map'.
\medskip

\noindent
{\bf Theorem A.4.}
(Multi-jet \emph{transversality} theorem of Thom,
 see \cite[section~I.2]{AVG})\\
Let $M,N$ be compact smooth manifolds,
 $\La\subset J^k_{[l]}(M,N)$ be a stratified space.
Given a smooth map $\xi:M\to N$, there is
 a smooth map $\eta:M\to N$ such that
\smallskip

$\bu$
 the map $\eta$ is arbitrarily close to $\xi$
  with respect to the Whitney topology;
\smallskip

$\bu$
 the $l$-tuple $k$-jet
 $j^k_{[l]}\eta:M^l\to J^k_{[l]}(M,N)$
 is transversal to $\La\subset J^k_{[l]}(M,N)$.
\qed
\medskip

\noindent
\emph{Proof of Claim~4.3}.
(i)
For any critical point of $\pr:S\to\R$, fix local coordinates $(x,y)\in S$
 such that the derivatives $\pr_x=\pr_y=0$.
The closures of the subspaces $\overline{\Si_{++}\cup\Si_{+-}\cup\Si_{--}}$ 
 and $\bar\Si_2$ from Definition~2.5 can be mapped onto the subspaces 
 of the finite-dimensional spaces $J_{[2]}^1(S,\R)$ and $J_{[1]}^2(S,\R)$
 given by the equations $\pr(x_1,y_1)=\pr(x_2,y_2)$ and 
 $\pr_{xx}\pr_{yy}-\pr_{xy}^2=0$, respectively.
The resulting subspaces of jets have codimension~1 
 as preimages of 0 under smooth functions, e.g. the image of $\bar\Si_2$ 
 in $J_{[1]}^2(S,\R)$ is $(\pr_{xx}\pr_{yy}-\pr_{xy}^2)^{-1}(0)$.
Hence the closures $\overline{\Si_{++}\cup\Si_{+-}\cup\Si_{--}}$ and
 $\bar\Si_2$ have codimension~1 in the space $\CS$ of 2-links.
\smallskip

\noindent
(ii)
If a 2-link is not generic, then either some critical points of 
 the projection $\pr:S\to\R$ are degenerate or have the same value.
The singularities of Definition~2.5 are all multi local codimension~1
 singularities of smooth functions $\R^2\to\R$, see \cite{AVG}. 
\smallskip

\noindent
(iii)
By Theorem A.4 any smooth isotopy of 2-links is a path in $\CS$
 and can be made transversal to the singular subspace $\bar\Si$, 
 which has codimension~1 by (i), hence the new path 
 will contain only finitely many isolated singularities of Definition~2.5.
\smallskip

\noindent
(iv)
The normal form of an $A_2$-singularity of a function $\R^2\to\R$
 is $\pr(x,y)=x^2-y^3$, i.e. the projection $\pr:S\to\R$ has 
 the form above in suitable local coordinates $(x,y)\in S$.
A 2-link $S$, its cross-sections around the singularity and 
 the graph of $y^3$ look like the middle pictures of Fig.~18.
The versal deformation of an $A_2$-singularity is $\pr(x,y;\ep)=x^2-y^3+\ep y$ 
 \cite{AVG}, i.e. any smooth deformation of $\pr(x,y)$ can be expressed as
 $f_1(x,y;\ep)\cdot\pr(f_2(x,y;\ep),f_3(x,y;\ep);f_4(\ep))$, 
 where $f_1,f_2,f_3,f_4$ are smooth, $f_1(0,0;0)\neq 0$, 
 $f_2(x,y;0)\equiv x$, $f_3(x,y;0)\equiv y$ and $f_4(0)=0$.
%\medskip

\begin{figure}[!h]
\includegraphics[scale=1.0]{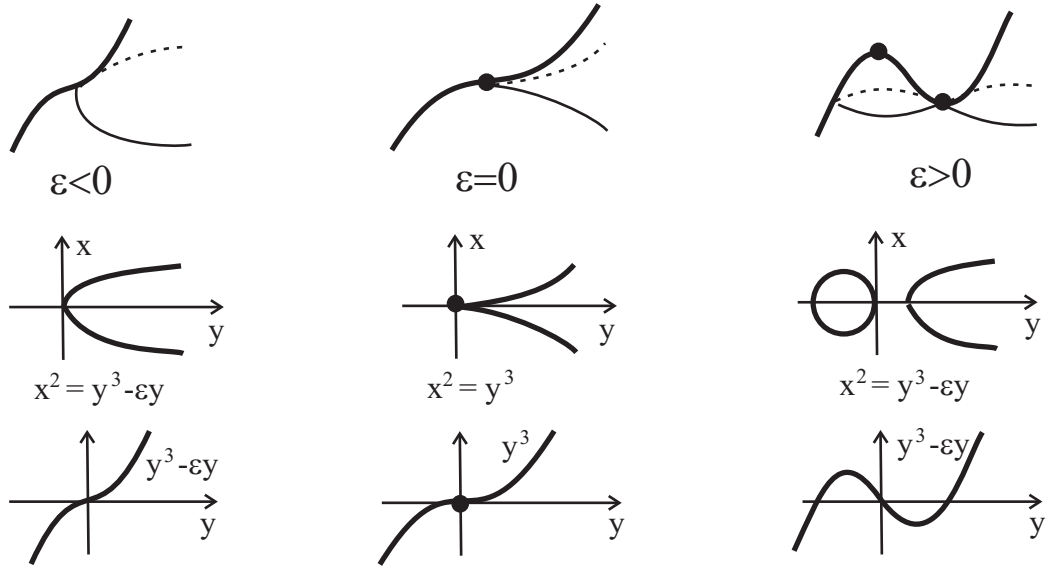}
\caption{Transformation of a 2-link near an $A_2$-singularity}
\end{figure}
\smallskip

For $\ep<0$, a 2-link $S$, its cross-sections around the singularity 
 and the graph of $y^3-\ep y$ look like the left pictures of Fig.~18.
For $\ep>0$, a 2-link $S$, its cross-sections around the singularity 
 and the graph of $y^3-\ep y$ look like the right pictures of Fig.~18.
For instance, 2-links for $\ep>0$ have a non-degenerate saddle at
 $x=0$, $y=\sqrt{\ep/3}$ and a local extremum at $x=0$, $y=-\sqrt{\ep/3}$.
\qed
\medskip

%======================================================


\begin{thebibliography}{References}

\bibitem{AVG}
\emph{V.~Arnold, A.~Varchenko, S.~Gusein-Zade},
Singularities of Differentiable Maps, Moscow, 1982.

\bibitem{CS}
\emph{S.~Carter, M.~Saito}, 
Knotted Surfaces and Their Diagrams,
Math. Surveys and Monographs, v.~55,
Providence, R.I. : Amer. Math. Soc., 1998.

\bibitem{FK}
\emph{T.~Fiedler, V.~Kurlin},
A one-parameter approach to links in solid torus,
math.GT/0606381. 

\bibitem{FM}
\emph{R.~Fox, J.~Milnor},
Singularities of 2-spheres in 4-space and cobordism of knots,
Osaka J. Math. {\bf 3} (1966), 257--267.

\bibitem{Kau}
\emph{L.~Kauffman},
Invariants of Graphs in Three-Space,
Trans. Amer. Math. Soc., {\bf 311} (1989), no.~2, 697--710.

\bibitem{Kaw}
\emph{A.~Kawauchi},
A survey of knot theory, Birkhauser Verlag (1996). 

\bibitem{KSS}
\emph{A.~Kawauchi, T.~T.~Shibuya, S.~Suzuki},
Description of surfaces in 4-space.I.Normal forms,
Math. Sem. Notes Kobe Univ., {\bf 10} (1982), 75--125.

\bibitem{KL}
\emph{C.~Kearton, W.~B.~R.~Lickorish},
Piecewise linear critical levels and collapsing,
Trans. Amer. Math. Soc., {\bf 170} (1972), 415--424.

\bibitem{Kur}
\emph{V.~Kurlin},
Three-page encoding and complexity theory for spatial graphs,
J. Knot Theory Ramifications {\bf 16} (2007), no.~1, 59--102.

\bibitem{KV}
\emph{V.~Kurlin, V.~Vershinin},
Three-page embeddings of singular knots,
Functional Anal. Appl. {\bf 35} (2004), no.~1, 14--27.

\bibitem{Swe}
\emph{F.~Swenton},
On a calculus for 2-knots and surfaces in 4-space,
J. Knot Theory Ramif. {\bf 10} (2001), no. 8, 1133--1141.

\bibitem{Yos}
\emph{K.~Yoshikawa},
An enumeration of surfaces in 4-space,
Osaka J. Math. {\bf 31} (1994), 497--522.

\end{thebibliography}
\end{document}